\documentclass[10pt]{article}
\usepackage{latexsym,amssymb,amsmath,amsfonts,amsthm,mathrsfs}

\def\cvd{\hfill$\Box$}
\def\vuoto{\varnothing}

\begin{document}

\begin{center}
\large \textbf{A remark on $\{2,3\}$-groups with no elements of order six} 
\normalsize \vspace{5mm}

\textsc{E. Jabara}

\end{center}

\vspace{5mm}

\begin{quote} {\bf Abstract:} {\it We describe $\{2,3\}$-groups in which the order of a 
product of any two elements of orders at most $4$ does not 
exceed $9$ and the centralizer of every involution is a 
locally cyclic $2$-subgroup. In particular, we will prove 
that these groups are locally finite.} \end{quote}

\vspace{5mm}

\textbf{\S 1. Introduction}

\vspace{3mm}

The purpose of this paper is to refine some results obtained 
by Lytkina and Mazurov in  \cite{LM} on $\{2,3\}$-groups 
with no elements of order $6$.

Let $G$  be a periodic group, we denote by $\pi(G)$
the set of primes that divide the order of some element 
of $G$ and $\Omega(G)=\langle \, g \mid g^{p}=1 
\mbox{ for some } p \in \pi(G) \,\rangle$. 
By $C_{n}$ we  denote the cyclic group of order $n$, 
by $D_{2n}$ the dihedral group of order $2n$ and by 
$S_{3} \simeq D_{6}$ the symmetric group on 
three symbols. 

We say that a group $G$ satisfies the hypotesis $({\bf H})$ 
if the following two conditions hold:

\begin{itemize}
\item[(H.1)] $\pi(G)=\{2,3\}$,
\item[(H.2)] the product of any two elements of $G$, whose orders 
are at most $4$, does not exceed $9$.
\end{itemize}

\vspace{2mm}

Our goal is to prove the following

\vspace{3mm}

\textsc{Theorem 1.}{\it Let $G$ be a group satisfying the hypothesis} 
({\bf H}) {\it and suppose that the centralizer of every involution of
$G$ is a locally cyclic $2$-group. Then $G=O_{3}(G) \rtimes C$
with $O_{3}(G)$ an abelian subgroup of exponent at most $9$ and
$C$ a locally cyclic $2$-group acting regulary on $O_{3}(G)$.
In particular $G$ is locally finite.}

\vspace{3mm} 

Thanks to the previous result it is possible to restate in a more 
satisfactory form the main result of \cite{LM}.

\vspace{3mm} 

\textsc{Theorem 2.} (\cite{LM}) {\it Let $G$ be a group satisfying the hypothesis} 
({\bf H}) {\it and with no elements of order $6$. Then one of the following 
is true: }
\begin{itemize}
\item[(1)] {\it $G=O_{3}(G) \rtimes T$,  where  $O_{3}(G)$ is abelian of 
exponent at most $9$ and $T$ either is a cyclic or locally cyclic $2$-group 
or is a quaternion group of order $8$ or $16$. In this case $G$ is locally
finite.}
\item[(2)] {\it $G=O_{2}(G) \rtimes R$, where $O_{2}(G)$ 
is  nilpotent of nilpotency class at most two and exponent at most $8$
and $R$ is a $3$-group with a unique subgroup of order $3$, 
acting freely on $O_{2}(G)$.}
\item[(3)] {\it $G=O_{2}(G) \rtimes D$, where $D$ is a dihedral group 
of order  $6$ or $18$ and $O_{2}(G)$ has exponent at most $4$ and 
nilpotency class at most $2$. In this case $G$ is locally finite.}
\end{itemize}

\vspace{3mm}

We want to state here two further consequences of Theorem 2

\vspace{3mm}

\textsc{Corollary 1.} {\it Let $G$ be a group satisfying the 
hypothesis} ({\bf H}) {\it and with no elements of order $6$, 
then $\Omega(G)$ is a locally finite 
group of exponent at most $72$.
Moreover $G/\Omega(G)$ is or a $3$-group or a $2$-group (cyclic or  
locally cyclic or isomorphic to $C_{2} \times C_{2}$ or to $D_{8}$) 
or $G/\Omega(G) \simeq S_{3}$.} 

\vspace{3mm}

\textsc{Corollary 2.} {\it Let $G$ be a group satisfying the 
hypothesis} ({\bf H}) {\it and with no elements of order $6$. 
If $G$ does not contain elements of 
order $27$, then $G$ is locally finite.}

\vspace{3mm}

\textbf{\S 2. Proof of the Theorem 1.}

\vspace{3mm}

In this section with $G$ we denote a group satysfing the 
hypotheses of the Theorem 1 and we define 
\[ \Gamma_{n}(G) =\{ g \in G \mid 
\mbox{$g$ has order exactly $n$} \}.\]
We will use freely the fact that if $K$ is a locally finite subgroup 
of $G$, then $K$ is a Frobenius group or a $2$-Frobenius group (that is 
$K=O_{3,2,3}(K)$ with $O_{3,2}(K)$ and $K/O_{3}(K)$ Frobenius groups).

To prove some lemmas we will use computations in \texttt{GAP} \cite{GAP}.

\vspace{3mm}

\textsc{Lemma 1.} {\it Let $a,b \in \Gamma_{2}(G)$, then $(ab)^{9}=1$.
In particular $G$ contains a single conjugacy class of involutions 
and $\Gamma_{2}(G)=a^{G}=a^{\Gamma_{2}(G)}$.}

\vspace{1mm}

\textsc{Proof.} Let $a,b$ be two involutions in $G$. Then $\langle
a,b \rangle \simeq D_{2n}$ and, since $G$ is $C_{2} \times C_{2}$-free,
then $n \in \{1,3,9\}$. The second part of the assertion follows from 
the fact that if $a,b \in \Gamma_{2}(D_{2n})$ with $n$ odd, then there is
$c \in \Gamma_{2}(D_{2n})$ such that $a=b^{c}$. \cvd

\vspace{3mm}

\textsc{Lemma 2.} {\it Let $a \in \Gamma_{2}(G)$, 
$x \in \Gamma_{3}(G)$
and $t \in \Gamma_{4}(G)$ with $x^{a}=x^{-1}$
and $t^{2}=a$. Then 
$\langle t,b \rangle \simeq (C_{3} \times C_{3}) \rtimes C_{4}$.}

\vspace{1mm}

\textsc{Proof.} Consider the groups
\[ H(i,j)=\langle t,x \mid \, t^{4}, \;\; x^{3}, \;\;
(t^{2}x)^{2}, \;\; (tx)^{i}, \;\; [t,x]^{j} \, \rangle \] 
($i,j \in \{8,9\}$). A computation
with \texttt{GAP} shows that $H(8,8) \simeq C_{4}$,
$H(8,9) \simeq (C_{3} \times C_{3}) \rtimes C_{4}$,
$H(9,8) \simeq PSL(2,19)$ and $H(9,9)=1$. \cvd

\vspace{3mm}

\textsc{Lemma 3.} {\it Let $x \in \Gamma_{3}(G)$ and 
$a \in \Gamma_{2}(G)$ with $x^{a}=x^{-1}$. Then $C_{G}(x)$ 
is an abelian $3$-group of exponent at most $9$ where $a$ 
acts as the inversion.}

\vspace{1mm}

\textsc{Proof.} Clearly $C=C_{G}(x)$ is an $a$-invariant $3$-subgroup of
$G$. Since in $K=C \rtimes \langle a \rangle$ is 
$C_{K}(a)=\langle a \rangle$, by [Ro] 14.3.7 we deduce that
$K$ is locally finite. From this fact the claim follows. \cvd

\vspace{3mm}

\textsc{Lemma 4.} {\it Let $a \in \Gamma_{2}(G)$, $x \in \Gamma_{9}(G)$ 
with $x^{a}=x^{-1}$ and let $t \in \Gamma_{4}(G)$ with
$t^{2}=a$. Then $\langle t,x \rangle$ is a Frobenius group isomorphic to 
a quotient of  $(C_{9} \times C_{9}) \rtimes C_{4}$. 
In particular $[x,x^{t}]=1$.}

\vspace{1mm}

\textsc{Proof.} Write $y=x^{t}$. By Lemma 2 $\langle x^{3}, t \rangle \simeq
(C_{3} \times C_{3}) \rtimes C_{4}$ and hence $[x^{3}, y^{3}]=1$. 
In particular $y^{3},x \in C_{G}(x^{3})$ and, by Lemma 3, $[x,y^{3}]=1$;
similarly $[x^{3},y]=1$ and hence $Z=\langle x^{3}, y^{3} \rangle
\leq Z(\langle x,y \rangle)$. Since $\langle x,y \rangle \leq C_{G}(Z) \leq
C_{G}(x^{3})$, by Lemma 3 we conclude that $\langle x,y \rangle$ is abelian.
Let \[ H(i,j)=\langle t,x \mid \, t^{4}, \;\; x^{9}, \;\;
(t^{2}x)^{2}, \;\; [x,x^{t}],\;\;  (tx)^{i}, \;\; [t,x]^{j} \, \rangle \] 
($i,j \in \{8,9\}$). A computation
with \texttt{GAP} shows that $H(8,8) \simeq C_{4}$,
$H(8,9) \simeq (C_{9} \times C_{9}) \rtimes C_{4}$ and
$H(9,8)=1=H(9,9)$. \cvd

\vspace{3mm}

We fix an element $t \in \Gamma_{4}(G)$ and we define 
\[ \Theta=t^{G}=\{t^{g} \mid g \in G \} \;\;\;\;\; \mbox{ and } 
\;\;\;\;\; \Theta^{-}=\{u^{-1} \mid u \in \Theta\}. \]

\vspace{3mm}

\textsc{Lemma 5.} {\it $\Theta \cap \Theta^{-}=\vuoto$ 
and $\Theta \cup \Theta^{-}=\Gamma_{4}(G)$.}

\vspace{1mm}

\textsc{Proof.} If $u \in \Theta \cap \Theta^{-}$, then there is $g \in G$
with $u^{g}=u^{-1}$. Clearly $g$ cannot have odd order and hence $g$ is a 
$2$-element of $G$. If $g \in \Gamma_{2}(G)$, then $\langle u^{2}, g \rangle
\simeq C_{2} \times C_{2}$, a contradiction. If $g \in \Gamma_{2^{n}}(G)$
with $n \geq 2$, then $\langle u, g \rangle \leq C_{G}(g^{2^{n-1}})$, a contradiction
since, by hypothesis, the centralizer of an involution in $G$ is a locally
cyclic $2$-group. 

Let $\Theta =t^{G}$ and let $u \in \Gamma_{4}(G)$ with $u^{2} \not =t^{2}$,
then $x=t^{2}u^{2} \in \Gamma_{3}(G) \cup \Gamma_{9}(G)$ and, by Lemma 4, 
$\langle t, x \rangle$ is a Frobenius group with complement 
$\langle t \rangle$ and kernel $N=\langle x, x^{t} \rangle$. 
Let $y \in N$ be such that $(t^{2})^{y}=u^{2}$.
Then $\langle t^{y} \rangle \leq C_{G}(u^{2})$ and hence $t^{y}=u$ or $t^{y}=u^{-1}$, 
since $C_{G}(u^{2})$ is locally cyclic. \cvd

\vspace{3mm}

\textsc{Lemma 6.} {\it Let $t,u \in \Theta$, then
$(tu)^{2}=1=(tu^{-1})^{9}$ and $(tu^{2})^{4}=1$.}

\vspace{1mm}

\textsc{Proof.} Let $x=t^{2}u^{2}$, then $F=\langle x,t \rangle$ is a
Frobenius group and, arguing as in the proof od Lemma 5, $u \in F$.
Since $\langle t \rangle$ and $\langle u \rangle$ are two complements
of $F$ and $x^{9}=1$, the claim follows. \cvd

\vspace{3mm}

\textsc{Lemma 7.} {\it Let $a,b,c \in \Gamma_{2}(G)$,
then $abc \in \Gamma_{2}(G)$.}

\vspace{1mm}

\textsc{Proof.}  Let $t_{1},t_{2},t_{3} \in \Theta$ be 
such that $t_{1}^{2}=a$, $t_{2}^{2}=b$
and $t_{3}^{2}=c$. By Lemma 6 in the group 
$K=\langle t_{1}, t_{2}, t_{3}  \rangle$ the following 
set of relations holds
\[ \mathcal{R}(t_{1},t_{2},t_{3})
=\big \{ (t_{\ell}t_{m}^{t_{n}})^{2}=(t_{\ell}^{-1}t_{m}^{t_{n}})^{9}=
(t_{\ell}^{2}t_{m}^{t_{n}})^{4}=1 \; \mid \; \ell, m, n \in \{1,2,3 \} \big\}. \]
A computation with \texttt{GAP} shows that 
$\langle \, t_{1}', t_{2}', t_{3}' \, \mid \mathcal{R}(t_{1}',t_{2}',t_{3}') \, 
 \,\rangle$ 
is isomorphic to the Frobenius group $(C_{9} \times C_{9} \times C_{9} \times C_{9}) 
\rtimes C_{4}$. Therefore $\langle t_{1}, t_{2}, t_{3} \rangle \leq G$ is a 
Frobenius group with complements $\langle t_{1} \rangle$, $\langle t_{2} \rangle$ 
and $\langle t_{3} \rangle$ and the lemma follows. \cvd

\vspace{3mm}

\textsc{Lemma 8.} {\it Let $R=\langle \Gamma_{2}(G) \rangle$, then
$R=O_{3}(R) \rtimes \langle a \rangle$ for every $a \in \Gamma_{2}(G)$.
Moreover $R$ is a Frobenius group and $O_{3}(R)$ is 
abelian of exponent at most $9$.}

\vspace{1mm}

\textsc{Proof.} Let $a_{1},a_{2}, \ldots, a_{n} \in \Gamma_{2}(G)$, 
then, by Lemma 7, $(a_{1}\cdot a_{2} \cdot \ldots \cdot a_{n})^{2}=1$ 
if $n$ is odd and $(a_{1}\cdot a_{2} \cdot \ldots \cdot a_{n})^{9}=1$ 
if $n$ is even.
Let $S$ be the subset of $G$ of the products of an even number of
elements of $\Gamma_{2}(G)$. By Lemma 7, if $s \in S$ and 
$a \in \Gamma_{2}(G)$, then $s^{a}=s^{-1}$ and since $s^{9}=1$, $S$ is
an abelian $3$-subgroup of $G$ of exponent at most $9$. Since 
$\Gamma_{2}(G)=a^{G}=\{a^{b} \mid b \in \Gamma_{2}(G) \}=\{a^{ab} \mid b \in 
\Gamma_{2}(G) \}$ and $ab \in S$ we can conclude that $R=S\langle a \rangle$
and $S=O_{3}(R)$. \cvd

\vspace{3mm}

\textsc{Proof of  Theorem 1.} Let $R$ and $S$ be the subgroups of $G$
defined in the proof of Lemma 8.
Since $R$ is normal in $G$ and $S$ is characteristic in 
$R$, the subgroup $S$ is normal in $G$. We have that $G/S$ contains a 
unique involution, hence $G/S$ is a locally cyclic $2$-group and 
$S=O_{3}(G)$. \cvd

\vspace{3mm}

\textbf{\S 3. Proof of the other results.}

\vspace{3mm}

\textsc{Proof of  Theorem 2.} Let $G$ be a group satisfying the 
hypothesis ({\bf H}) and suppose that $G$ does not contain elements 
of order $6$. By Theorem 1, $G$ cannot contains {\em forbidden} 
subgroups as defined in \cite{LM} (a forbidden subgroup $Y$ is a subgroup 
of $G$ generated by an involution and an element of order $3$ and 
such that every maximal $2$-subgroup of $Y$ is locally cyclic). 
The main theorem of \cite{LM} tells us that $G$ must be as described 
in the statement of Theorem 2.

The bounds to the exponents of the subgroups of $G$ in claims 
(1), (2) and (3) are 
trivial consequences of the hypothesis (H.2). \cvd

\vspace{3mm}

\textsc{Proof of  Corollary 1.} A straightforward consequence 
of Theorem 2. \cvd

\vspace{3mm}

\textsc{Proof of  Corollary 2.} Let $G$ be a group satysfing the 
hypothesis of the Corollary 2. The only case in which
$G$ can be not locally finite is when $G$ is as described in the 
claim (2) of the Theorem 2.
In this case, if $\Gamma_{27}(G)=\vuoto$, then $G/\Omega(G)$ is 
a group of exponent $3$ and hence locally finite (see \cite{Ro} 14.2.3). 
By Corollary 1 also $\Omega(G)$ is locally finite, 
therefore, by \cite{Ro} 14.3.1, $G$ itself is locally finite. \cvd

\vspace{5mm}

\textsc{Enrico Jabara}

\textsc{DFBBC Universit\`{a} di Ca' Foscari}

\textsc{Dorsoduro 3484/D -- 30123 Venezia -- ITALY}

{\it Email address:} \texttt{jabara@unive.it}

\end{document}